\documentclass[12pt,leqno]{article}
\usepackage{amsfonts}
\usepackage{amssymb}
\usepackage{amsthm}
\usepackage{amsmath}
\usepackage{graphicx}

\newtheorem{problem}{Problem}[section]
\newtheorem{theo}[problem]{Theorem}

\newtheorem{cor}[problem]{Corollary}

\begin{document}

\date{October, 2007}
 \title{{\Large \bf Oriented matroids\\ and Ky Fan's theorem
 }}
\author{Rade  T. \v Zivaljevi\' c\\
Mathematical Institute SANU, Belgrade}

\maketitle
\begin{abstract}\noindent
L.~Lov\' asz has shown in \cite{Lov80} that Sperner's
combinatorial lemma admits a generalization involving a {\em
matroid} defined on the set of vertices of the associated
triangulation. We prove that Ky Fan's theorem admits an {\em
oriented matroid} generalization of similar nature
(Theorem~\ref{thm:fanki-matroid}). Classical Ky Fan's theorem is
obtained as a corollary if the underlying oriented matroid is
chosen to be the alternating matroid $C^{m,r}$.
\end{abstract}

\section{Introduction}\label{sec:intro}

The following extension of Sperner's combinatorial lemma was
proposed by L\' aszlo Lov\' asz in \cite{Lov80}.

\begin{theo}{{\em (\cite{Lov80})}}\label{thm:Lovasz}
Let $K$ be a simplicial complex which triangulates a
$d$-dimensional manifold. Suppose that a matroid $\mathcal{N}$ of
rank $d+1$ is defined on the set ${\rm vert}(K)$ of vertices of
$K$. If $K$ has a simplex whose vertices form a basis of the
matroid $\mathcal{N}$, than there exist at least two such
simplices.
\end{theo}

\noindent Homological nature of this result was subsequently
emphasized by Bernt Lindstr\" om who demonstrated that it suffices
to assume that the complex $K$ (in Theorem~\ref{thm:Lovasz})
supports a $d$-cycle.

\begin{theo}{\em (\cite{Lind81})}
Let $C=\sum_{i\in I}\alpha_i\sigma_i$ be a $d$-cycle
in a complex $K$. Assume that the vertices of simplices in $C$ are
labelled by elements of a matroid $\mathcal{N}$ of rank $(d+1)$.
If some simplex is labelled by the elements of a base of
$\mathcal{N}$, then there are at least two simplices in $C$ with
this property.
\end{theo}

\noindent As a corollary of Theorem~\ref{thm:Lovasz} Lov\' asz
deduced the following result which reduces to the classical
version of Sperner's lemma if $\mathcal{N}$ is the matroid such
that $S\subset {\rm vert}(K)$ is an independent set if and only if
its elements are labelled (colored) by different labels.

\begin{cor}
Suppose that $K$ is a simplicial subdivision of a $d$-dimensional
simplex $\sigma$. Suppose that a matroid $\mathcal{N}$ of rank
$d+1$ is defined on the set ${\rm vert}(K)$ of vertices of $K$ and
let ${\rm vert}(\sigma)$, the set of vertices of $\sigma$, be a
basis of $\mathcal{N}$. Assume that for each face $F$ of $\sigma$,
the set $F\cap {\rm vert}(K)$ is in the flat of the matroid
$\mathcal{N}$ spanned by ${\rm vert}(F)\subset {\rm
vert}(\sigma)$. Then $K$ has a simplex whose vertices form a basis
of $\mathcal{N}$.
\end{cor}

\noindent Well known $\mathbb{Z}_2$-counterparts of Sperner's lema
are Tucker's lemma \cite{Tuc45}, \cite{Mat03} and its
generalization due to Ky Fan \cite{Fan52}.

\begin{theo}{\em (\cite{Fan52})}\label{thm:fan52}
Suppose that $K$ is a $\mathbb{Z}_2$-invariant triangulation of
the sphere $S^n$. Let $\lozenge^m:=\mbox{\rm conv}\{\pm
e_1,\ldots,\pm e_m\}$ be the $m$-dimensional crosspolytope and
$\partial(\lozenge^m)\cong S^{m-1}$ its boundary with inherited
{\em (}$\mathbb{Z}_2$-invariant{\em )} triangulation. If $f :
K\rightarrow \partial(\lozenge^m)$ is a simplicial,
$\mathbb{Z}_2$-equivariant map, then $n<m$ and
$$
\sum_{1\leq k_1<k_2<\ldots<k_{n+1}\leq
m}\alpha(k_1,-k_2,k_3,-k_4,\ldots,(-1)^n k_{n+1}) \cong 1 \quad
(\mbox{\em mod}\, 2)
$$
where $\alpha(j_1, j_2,\ldots, j_{n+1})$ is the number of
$n$-simplices in $K$ mapped to the simplex spanned by vectors
$e_{j_1}, e_{j_2},\ldots, e_{j_{n+1}}$ and by definition $e_{-j}:=
- e_j$.
\end{theo}

\noindent A natural question is whether there exists a counterpart
of Lov\' asz' theorem  which extends Ky Fan's theorem (Tucker's
theorem) in the manner Theorem~\ref{thm:Lovasz} extends Sperner's
lemma. Our objective is to prove such a result,
Theorem~\ref{thm:fanki-matroid}. It is not a surprise that
oriented matroids appear in this extension and play a role similar
to the role matroids play in Theorem~\ref{thm:Lovasz}.

\section{Oriented matroids in a nutshell}

Oriented matroids provide combinatorial models for important
geometric objects, configurations and structures including the
following:

\begin{itemize}
 \item linear subspaces $L\subset \mathbb{R}^n$, \vspace{-0.2cm}
 \item configurations of points (vectors) in $\mathbb{R}^n$, \vspace{-0.2cm}
 \item matrices, \vspace{-0.2cm}
 \item directed graphs, \vspace{-0.2cm}
 \item convex polytopes, \vspace{-0.2cm}
 \item linear programs, \vspace{-0.2cm}
 \item hyperplane arrangements etc.
\end{itemize}

\noindent Although they appear in many incarnations and disguises,
oriented matroids always provide essentially the same amount of
information about the object they discretize (cryptomorphism). The
reader is referred to \cite{Zie95} (Section~6) for a quick
introduction and initial motivation and to \cite{B+} for
illuminating orientation sessions (Sections~1 and 2) and thorough
treatment of the general theory with many interesting
applications. More recent reference \cite{RG-Z} offers both an
outline of the theory and a guide to the papers published after
the appearance of \cite{B+}.

\subsection{$\lozenge^m$-oriented matroid of a linear subspace}
\label{sec:lozenge}

Let us briefly review how an oriented matroid can be associated to
a $d$-dimensional linear subspace  $L\subset\mathbb{R}^m$. Let
$\lozenge^m:=\mbox{\rm conv}\{\pm e_1,\ldots,\pm e_m\}$ be the
crosspolytope in $\mathbb{R}^m$ and let $P(\lozenge^m)$ be the
associated face poset. Define $\lozenge^m$-oriented matroid
$L_{\lozenge^m}$ of $L$ as the collection (poset) $L_{\lozenge^m}
= \{F\in P(\lozenge^m) \mid L\cap F\neq\emptyset\}$ of all faces
of the crosspolytope intersected by $L$. The face poset
$P(\lozenge^m)$ is isomorphic to the poset $\{0,+,-\}^m$ of all
sign-vectors of length $m$. This isomorphism associates to each
face $F\in P(\lozenge^m)$ the sign vector $sign(F) :=
sign(v)\in\{0,+,-\}^m$ for some (any) $v\in relint(F)$. It follows
that $L_{\lozenge^m}$ is essentially the set $\mathcal{V}^\ast$ of
all {\em covectors} of an oriented matroid $\mathcal{M} =
\mathcal{M}(L)$ which captures the combinatorial information about
how the subspace $L$ is placed in the ambient space
$\mathbb{R}^m$. Similarly, the set $\mathcal{C}^\ast$ of all {\em
cocircuits} of $\mathcal{M}(L)$ can be described as the set of all
$\subseteq$-minimal elements in $L_{\lozenge^m}$ or alternatively
as the collection of faces $F\in P(\lozenge^m)$ such that the
subspace $L$ intersects $relint(F)$ in a single point.

\subsection{Topological representation theorems}

The original {\em Topological Representation Theorem} for oriented
matroids was proved by Folkman and Lawrence \cite{B+}. The
following strengthening of this result, due to Brylavski and
Ziegler \cite{BryZie}, provides a simultaneous representation for
both the oriented matroid $\mathcal{M}$ and its dual
$\mathcal{M}^\ast$.

\begin{theo}{\rm (\cite{BryZie})}\label{thm:dual}
For each oriented matroid  $\mathcal{M}$
of rank $r$ on $\{1,\ldots,m\}$ there exists a (signed)
pseudosphere arrangement $\mathcal{A}=(S_i)_{1\leq i\leq 2m}$ in
$S^{m-1}$ such that:
\begin{enumerate}
 \item[{\rm (1)}] $S_i=\{x\in S^{m-1}\mid
 x_i=0\}$ for $1\leq i\leq m$. \vspace{-0.3cm}
 \item[{\rm (2)}] The $(r-1)$-subsphere $S_B:=S_{m+r+1}\cap\ldots \cap
 S_{2m}$ and $(m-r-1)$-subsphere $S_A:=S_{m+1}\cap\ldots \cap
 S_{m+r}$ are disjoint.\vspace{-0.3cm}
 \item[{\rm (3)}] The arrangement $(S_i\cap S_B)_{1\leq i\leq m}$
 is a topological representation of the oriented matroid $\mathcal{M}$
 in $S_B$.\vspace{-0.3cm}
 \item[{\rm (4)}] The arrangement $(S_j\cap S_A)_{1\leq j\leq m}$
 is a topological representation of the oriented matroid
 $\mathcal{M}^\ast$ in $S_A$.
\end{enumerate}
\end{theo}

\begin{cor}\label{cor:dual}
Let $\mathcal{M}$ be an oriented matroid of rank $r$ on
$[m]=\{1,2,\ldots,m\}$. Let $S_i,\, i=1,\ldots,m$ be the
coordinate sphere $S_i=\{x\in S^{m-1}\mid x_i=0\}$. Then there
exists a $(r-1)$-dimensional pseudosphere $S_B$ in $S^{m-1}$ such
that the arrangement $\{S_i\cap S_B\}_{i=1}^m$ is a topological
representation of $\mathcal{M}$ in $S_B$. Moreover $S_B$ can be
chosen to be centrally symmetric and transverse to each of the
spheres $S_I:=S_{i_1}\cap\ldots\cap S_{i_{r}}$ for any $r$-element
subset $$I=\{i_1,\ldots, i_{r}\}\subset [m].$$
\end{cor}

\medskip\noindent
{\bf Proof:} According to Section~5.2 in \cite{B+}
(Theorem~5.2.1), the condition that all pseudospheres are
centrally symmetric can be always satisfied. In particular all
pseudospheres in Theorem~\ref{thm:dual} can be assumed to have
this property. Also, the transversality condition from
Corollary~\ref{cor:dual} is ``built in'' the Topological
Representation Theorem for oriented matroids. \hfill$\square$

\section{Generalization of Ky Fan's theorem}

\begin{theo}\label{thm:fanki-matroid}
Suppose that $M$ is a connected, $n$-dimensional, triangulated
$\mathbb{Z}_2$-manifold. Moreover, it is assumed that the
involution $\nu : M\rightarrow M$ defining the
$\mathbb{Z}_2$-action is fixed-point-free. Given a positive
integer $m$, let
$$\lambda : {\rm vert}(M)\rightarrow \{\pm 1,\pm 2,\ldots,\pm m\}$$
be a labelling of the vertices of $M$ which satisfies the
conditions:
\begin{enumerate}
 \item[{\rm (a)}]  $\lambda(\nu(v))=-\lambda(v)$ for each $v\in {\rm
vert}(M)$,
 \item[{\rm (b)}] $\lambda(v)+\lambda(v')\neq 0$ for each $1$-simplex
$\tau=\{v,v'\}$ in $M$.
\end{enumerate}
Let $\mathcal{M}=\mathcal{M}([m],\mathcal{C}^\ast)$ be an oriented
matroid of rank $r=m-n$ on the set $E=[m]=\{1,\ldots,m\}$ where
$\mathcal{C}^\ast$ is the set of associated cocircuits. Moreover,
$\mathcal{M}$ is assumed to be {\em uniform} in the sense that all
cocircuits in $\mathcal{C}^\ast$ have the same cardinality
$n=m-r$. Let $w_1\in H^1_{\mathbb{Z}_2}(M, \mathbb{Z}_2)\cong
H^1(M/\mathbb{Z}_2, \mathbb{Z}_2)$ be the first Stiefel-Whitney
class of the $\mathbb{Z}_2$-manifold $M$ and let $w(M)=
w_1^n([M/\mathbb{Z}_2])$ be the associated Stiefel-Whitney number
where $[M/\mathbb{Z}_2]$ is the $\mathbb{Z}_2$-fundamental class
of ${M}/\mathbb{Z}_2$. Then,
\begin{equation}\label{eqn:mod2}
w(M) = \frac{1}{2}\sum_{\tau\in \mathcal{C}^\ast} \alpha(\tau) =
\sum_{[\tau]\in \mathcal{C}^\ast/\mathbb{Z}_2} \alpha(\tau) \qquad
{\rm (mod\,\, 2)}
\end{equation}
where $\alpha(\tau)$ is the number of $n$-simplices $\sigma\in M$
whose vertices receive labels from $\tau$, i.e.\ such that $\tau =
\{\lambda(v)\mid v\in {\rm vert}(\sigma)\}$.
\end{theo}

\medskip\noindent
{\bf Proof:} Let $\partial\lozenge^m$ be the boundary of the
crosspolytope $\lozenge^m={\rm conv}\{\pm e_1,\ldots,\pm e_m\}$.
As a $\mathbb{Z}_2$-space, $\partial\lozenge^m$ is isomorphic to
the join $[2]\ast\ldots\ast [2]$ of $m$ copies of $[2]=\{1,2\}$.
We may therefore see $\partial\lozenge^m$ as a subcomplex of (the
boundary of) an ``infinite crosspolytope''
$E\mathbb{Z}_2=[2]\ast\ldots\ast[2]\ast\ldots$, which is nothing
but the well known Milnor's model for the classifying space of the
group $\mathbb{Z}_2$.

Any labelling $\lambda : {\rm vert}(M)\rightarrow \{\pm
1,\ldots,\pm m\}$ can be in an unique way extended to a
$\mathbb{Z}_2$-equivariant map
$$
\Lambda : M\rightarrow \lozenge^m
$$
which is affine on each simplex $\sigma\in M$. Conditions (a) and
(b) on the labelling guarantee that ${\rm
Im}(\Lambda)\subset\partial\lozenge^m$ and that $\Lambda$ is a
simplicial $\mathbb{Z}_2$-map. Consequently, $\Lambda$ is
essentially a classifying map $\Lambda: M\rightarrow
\partial\lozenge^m\hookrightarrow E\mathbb{Z}_2$, that is the
unique (up to a $\mathbb{Z}_2$-homotopy)
$\mathbb{Z}_2$-equivariant map $\Lambda: M\rightarrow
E\mathbb{Z}_2$. Let $\xi =\Lambda/\mathbb{Z}_2 :
M/\mathbb{Z}_2\rightarrow
\partial\lozenge^m/\mathbb{Z}_2$ be the induced map. Let
$\gamma\in
H^1_{\mathbb{Z}_2}(\partial\lozenge^m/\mathbb{Z}_2;\mathbb{Z}_2)\cong
H^1(\partial\lozenge^m/\mathbb{Z}_2; \mathbb{Z}_2)$ be the first
Stiefel-Whitney class of the $\mathbb{Z}_2$-space
$\partial\lozenge^m$, or equivalently the first S-W-class of the
line bundle
$$
\mathbb{R}^1\rightarrow
\partial\lozenge^m\times_{\mathbb{Z}_2} \mathbb{R}^1\rightarrow
\partial\lozenge^m/\mathbb{Z}_2 .
$$
Let $S_B$ be the $(r-1)$-dimensional pseudosphere in
$S^{m-1}\subset \mathbb{R}^m$ which represents the oriented
matroid $\mathcal{M}$ in the sense of Corollary~\ref{cor:dual}.
The fundamental class $y = [S_B/\mathbb{Z}_2]$ of the quotient
projective space $S_B/\mathbb{Z}_2\subset RP^{m-1}$ represents the
generator of the group $H_{r-1}(RP^{m-1};\mathbb{Z}_2)\cong
\mathbb{Z}_2$. This follows from the fact that, according to the
definition of a (centrally symmetric) pseudosphere, there exists a
$\mathbb{Z}_2$-homeomorphisms $h : S^{m-1}\rightarrow S^{m-1}$
which maps the $(r-1)$-dimensional equator of $S^{m-1}$ to the
pseudosphere $S_B$. By the well-known formula, see e.g.\
\cite{Dold}, Section~VII.12, or \cite{Bredon}, Theorem~VI.5.2 part
(4), and standard (loc.\ cit.) properties of the products of
(co)homology classes and the Poincar\' e duality map $D$,
\begin{equation}
w(M)= w_1^n([M/\mathbb{Z}_2]) = w_1^n\cap [M/\mathbb{Z}_2]=
\xi_\ast (\xi^\ast(\gamma^n)\cap [M/\mathbb{Z}_2])=
\end{equation}
\begin{equation}
=\gamma^n\cap \xi_\ast[M/\mathbb{Z}_2] = D(\gamma^n)\bullet
\xi_\ast[M/\mathbb{Z}_2] = [S_B]\bullet \xi_\ast[M/\mathbb{Z}_2].
\end{equation}
For the completion of the proof it is sufficient to observe that
the intersection product $[S_B]\bullet \xi_\ast[M/\mathbb{Z}_2]$
is precisely the right hand side of the equation (\ref{eqn:mod2}).
Indeed, the intersection $\tau\cap S_B$ is transverse for each of
the simplices (cocircuits) $\tau\in\mathcal{C}^\ast$ and each of
them is counted with the multiplicity $\alpha(\tau)$. \hfill
$\square$

\begin{cor}\label{cor:sphere}
If $w_1^n$ is non-trivial, i.e.\ if
$w(M)=w_1^n([M/\mathbb{Z}_2])=1$, for example if $M=S^n$ is the
$n$-sphere, then
\begin{equation}\label{eqn:cor:sphere}
1 = \sum_{[\tau]\in \mathcal{C}^\ast/\mathbb{Z}_2} \alpha(\tau)
\qquad {\rm (mod\,\, 2)}
\end{equation}
where the summation is over the representatives of classes
$[\tau]=\{\tau, -\tau\}$ of antipodal simplices (cocircuits).
\end{cor}
\begin{cor}\label{cor:fanky}
The usual Ky Fan's theorem (Theorem~\ref{thm:fan52}) follows from
Theorem~\ref{thm:fanki-matroid} (Corollary~\ref{cor:sphere}) if
$\mathcal{M}$ is chosen as (the dual of) the alternating oriented
matroid $C^{m,r}$.
\end{cor}

\medskip\noindent
{\bf Proof:} The (dual of) the alternating matroid $C^{m,r}$ is
defined, \cite{B+} Section~9.4, as the oriented matroid of the
vector configuration $\mathcal{W}=\{v_1,\ldots, v_m\}$ where
$v_i:=(1,t_i,\ldots,t_i^r)$ are the points on the moment curve
corresponding to a sequence $0<t_1<\ldots< t_m$. The associated
set of cocircuits is, following the description given in
Section~\ref{sec:lozenge}, obtained as the $\lozenge^m$-oriented
matroid associated to the subspace $L:={\rm Im}(W)$ where $W$ is
the matrix
$$
W = [v_1^t,\ldots, v_m^t] = \left [
 \begin{array}{cccc}
 1 & 1 & \dots & 1\\
 t_1 & t_2 & \dots & t_m\\
 \vdots & \vdots & \ddots & \vdots\\
 t_1^{r-1} & t_2^{r-1}&\dots &t_m^{r-1}
 \end{array}
 \right]
$$
It is well known and not difficult to prove that the associated
cocircuits are precisely the alternating sequences (and their
antipodes) that appear in the formulation of
Theorem~\ref{thm:fan52}. \hfill$\square$

\section{Homological reformulation}

Theorem~\ref{thm:fanki-matroid} admits a reformulation which
emphasizes its homological nature. It is parallel to and, to some
extent, inspired by Lindstr\" om's extension of the result of
Lov\' asz. Theorem~\ref{thm:comb-formula} can be interpreted as a
combinatorial formula (involving oriented matroids) for a power of
the first Stiefel-Whitney characteristic cohomological class of
the $\mathbb{Z}_2$-complex $\partial\lozenge^m\cong S^{m-1}$.
Moreover, it points in the direction of a hypothetical
``homological representation theorem'' for oriented matroids.

As before, each cocircuit (covector) of an oriented matroid
$\mathcal{M}$ is identified with the corresponding face $\tau$ of
the cross-polytope $\lozenge^m$ (Section~\ref{sec:lozenge}).

\begin{theo}\label{thm:comb-formula}
Let $\mathcal{M}=\mathcal{M}([m],\mathcal{C}^\ast)$ be an uniform
oriented matroid of rank $r=m-n$ on the set $E=[m]=\{1,\ldots,m\}$
where $\mathcal{C}^\ast$ is the set of associated cocircuits. For
each cocircuit $\tau\in \mathcal{C}^\ast\subset P(\lozenge^m)$ let
$\hat\tau$ be the $\mathbb{Z}_2$-cochain dual to $\tau$ where
$\hat\tau(\tau)=1$ and $\hat\tau(\theta)=0$ if $\theta\neq\tau$.
Then the $n$-dimensional, $\mathbb{Z}_2$-equivariant cochain

\begin{equation}\label{eqn:cochain}
C_{\mathcal{M}}:= \sum_{\tau\in \mathcal{C}^\ast} \hat\tau \in
C^n_{\mathbb{Z}_2}(\partial\lozenge^m)
\end{equation}
represents the Stiefel-Whitney characteristic class $$w_1^n\in
H^n_{\mathbb{Z}_2}(\partial\lozenge^m)\cong H^n(RP^{m-1}).$$
\end{theo}

\medskip\noindent
{\bf Proof:} The proof uses similar ideas as the proof of
Theorem~\ref{thm:fanki-matroid} so we omit the details. The key
observation is that the cochain $C_{\mathcal{M}}$ is a Poincar\' e
dual (on the level of chains) of the class $[S_B]$. For this it is
sufficient to check that $C_{\mathcal{M}}$ is a cochain which can
be deduced from the fact that the pseudosphere $S_B$ is transverse
to all the simplices in the triangulation of $\partial\lozenge^m$.
\hfill$\square$

\medskip\noindent
{\bf Example:} Choose a rank $2$ oriented matroid $\mathcal{M}$
associated to a nonstretchable arrangement $\mathcal{A}$ of
pseudolines. For example let $\mathcal{A}$ be the arrangement of
nine pseudolines described by Ringel, see \cite{RG-Z}
Figure~6.1.2. In this case $m=9, r=3$ and $n=9-3=6$. According to
formula (\ref{eqn:cochain}) this yields a cochain representative
for the class $w_1^6$ which has a $6$-dimensional simplex (and its
antipode) for each intersection of two pseudolines in
$\mathcal{A}$.

\section{Concluding remarks}

The proof of Theorem~\ref{thm:fanki-matroid} appears to be new
already in the realizable case where we don't need the full power
of the Topological representation theorem. In particular this
approach yields a short and conceptual, albeit non-constructive,
proof of the classical Ky Fan's theorem.

\par It is natural to ask if the condition (Theorem~\ref{thm:fanki-matroid})
that $\mathcal{M}$ is a uniform oriented matroid can be relaxed.
Indeed, it would be desirable and hopefully not too difficult to
come up with analogues of the formula (\ref{eqn:cochain}) for the
case of non-uniform matroids. The critical step is to express the
intersection product $[S_B]\bullet \xi_\ast[M/\mathbb{Z}_2]$ in
terms of the underlying oriented matroid. One way around this
difficulty is to analyze perturbations of the associated
pseudosphere arrangement.

Aside from the conceptual interest, more general formulas could
lead to new inductive and constructive proofs based on standard
oriented matroid technique. This may prove useful in finding new
systematic ways of producing combinatorial proofs for
combinatorial statements which originally required topological
methods, cf.\ \cite{LoZiv06} \cite{Mat04} \cite{Meu06}
\cite{Zie02} for some of the more recent related developments.

Considering some recent advances \cite{Meu05} \cite{Zie02} in
understanding $Z_q$-analogues of Tucker's and Ky Fan's theorem, it
would be interesting to know if such analogs exist for
Theorem~\ref{thm:fanki-matroid}. This may involve a development of
an analogue (or replacement) for the concept of a
($\mathbb{Z}_q$-oriented) matroid, see \cite{Zie93} for a related
development.

Formula (\ref{eqn:mod2}) seems to indicate that, at least in
principle, all the formulas involving algebraic count of {\em
alternating simplices}, could be instances of more general
statements involving oriented matroids. For example it is
plausible that Sarkaria's ``Generalized Tucker-Ky Fan theorem''
\cite{Sar90} admits such a generalization.

Finally, in light of the fact that Ky Fan's theorem and its
consequences have found numerous applications in combinatorics and
discrete geometry, \cite{ST06} being one of the latest examples,
it remains to be seen if Theorem~\ref{thm:fanki-matroid} can be
used for a similar purpose.

\leftskip 0cm \hoffset 0cm

\end{document}